\documentclass[a4paper,11pt]{article}
\usepackage{amssymb,palatino}
\newtheorem{theorem}{Theorem}

\newtheorem{lemma}{Lemma}
\newtheorem{definition}{Definition}

\newcommand{\hhs}[1]{\hspace{#1mm}}

\newcommand{\sty}{\displaystyle}
\newcommand{\QED}{\begin{flushright} $\Box$ \end{flushright}}
\newcommand{\gauss}[2]{\left[ \begin{array}{c} #1\\ \hline #2 \end{array} \right]}
\newcommand{\gaussf}[2]{\left[ \begin{array}{c} #1\\ \hline #2 \end{array} \right]_4}
\newcommand{\jacobi}[2]{\left( \begin{array}{c} #1\\ \hline #2 \end{array} \right)}

\newcommand{\dede}[2]{\left( \left( \begin{array}{c} #1\\ \hline #2 \end{array}\right) \right)}
\newcommand{\deder}[1]{(( #1  ))}
\newcommand{\notdiv}{{ \hspace{1.1mm} | \hspace{-1.5mm}/\hhs{2}}}
\begin{document}

\title{The Rabin cryptosystem revisited}

\author{Michele Elia\thanks{Politecnico di Torino, Italy}, ~~
 Matteo Piva \thanks{Universit\'a di Trento, Italy}, ~~
 Davide Schipani \thanks{University of Zurich, Switzerland}
}

%\date{August 2010}
\maketitle

\thispagestyle{empty}

\begin{abstract}
\noindent
The Rabin public-key cryptosystem is revisited with a focus on the problem of identifying the encrypted message unambiguously for any pair of primes. In particular,
 a deterministic scheme using quartic reciprocity is described that 
 works for primes congruent $5$ modulo $8$, a case that was still open.
Both theoretical and practical solutions are presented. 
The Rabin signature is also reconsidered and a deterministic padding mechanism is proposed. 
\end{abstract}


\vspace{2mm}

% ********************************************************************
\vspace{8mm}

\section{Introduction}
In 1979, Michael Rabin \cite{rabin} suggested a variant of RSA with public-key
 exponent $2$, which he showed to be as secure as factoring. 
The encryption of a message $m \in \mathbb Z_N^{*}$ is 
$C = m^2  \bmod N $,  
 where  $N=p q$ is a product of two prime numbers,
and decryption is performed by solving the equation
\begin{equation}
  \label{main}
  x^2 = C  \bmod N ~~, 
\end{equation}
which has four roots; thus for complete decryption, further information is needed 
 to identify $m$ among these roots. 
More precisely, for a fully automatic (deterministic)
 decryption we need at least two more bits (computed at the encryption stage) 
 to identify $m$ without ambiguity. 
The advantages of using this exponent $2$, compared to larger exponents, are: i) a smaller
 computational burden, and ii) solving (\ref{main}) is equivalent to factoring $N$.
The disadvantages are: iii) computation, at the encryption stage, of the information required
 to identify the right root, and the delivery of this information to the decryption stage,
 and iv) vulnerability to chosen-plain text attack \cite{buchmann,menezes3,schneier,silverman}. 
Several naive choice methods base selection of the correct root
 on the message semantics, that is they retain the root that corresponds 
 to a message that looks most meaningful, or the root that contains a known string of bits. 
However, these methods are either unusable, for example when
 the message is a secret key, or are only probabilistic; in any
 case they affect the equivalence between breaking the Rabin scheme and factoring \cite{buchmann}. 
Nevertheless, for schemes using pairs of primes congruent $3$ modulo $4$ (Blum primes),
 Williams \cite{williams} proposed a root identification scheme based on the computation of a
 Jacobi symbol, using an additional parameter in the public key, and two additional bits in
 the encrypted message. 

The Rabin cryptosystem may also be used to create a signature by exploiting the inverse mapping:
 in order to sign $m$, the equation $x^2 =m \bmod N$ is solved and any of the four roots, say $S$,
 can be used to form the signed message $(m,S)$. However, if $x^2 =m \bmod N$ has no solution,
 the signature cannot be directly generated; to overcome this issue, a random pad $U$
 is used until $x^2 =mU \bmod N$ is solvable, and the signature is the triple
 $(m,U,S)$ \cite{seberry}. A verifier compares $S^2$ with $mU \bmod N$ and accepts the signature as valid
 when these two  numbers are equal. For an  
 application to electronic signature, an in-depth analysis on advantages/disadvantages can be found in \cite{bernstein}.

The next Section provides preliminary results concerning the solutions
 of the equation (\ref{main}) and the mathematics that will be needed. 
Section 3 describes in detail the Rabin scheme in the standard setting,
 where both prime factors of $N$ are congruent $3$ modulo $4$, 
 and proposes a new identification rule exploiting the Dedekind sums.
Section 4 addresses the identification problem for any pair of primes, featuring a deterministic scheme working with primes 
 congruent $5$ modulo $8$ based on quartic residues of Gaussian integers.
Section 5 considers a Rabin signature with deterministic padding.
Lastly, Section 6 draws some conclusions.
 
\section{ Preliminaries}
Let $N=pq$ be a product of two odd primes $p$ and $q$. Using the
 generalized Euclidean algorithm to compute the greatest common divisor between $p$ and $q$,
 two integer numbers, $\lambda_1,\lambda_2 \in \mathbb Z$, such that $\lambda_1 p+ \lambda_2 q = 1$, are
 efficiently computed. Thus, setting
 $\psi_1 =\lambda_2  q$ and $\psi_2 =\lambda_1 p$, so that $\psi_1+\psi_2=1$,
it is easily verified that $\psi_1$ and $\psi_2$ satisfy the relations  %"ortho-normality" 
\begin{equation}
\label{normal}
        \left\{ \begin{array}{l} 
             \psi_1 \psi_2 =0 \bmod N  \\
             \psi_1^2= \psi_1 \bmod N  \\
             \psi_2^2= \psi_2 \bmod N ~~. \\
               \end{array}       
              \right.           
\end{equation} 
and that $\psi_1= 1 \bmod p$, $\psi_1= 0 \bmod q $, and $\psi_2= 0 \bmod p$, $\psi_2= 1 \bmod q $.
According to the Chinese Remainder Theorem (CRT), using $\psi_1$ and $\psi_2$, every element $a$
 in $\mathbb Z_N$ can be represented as
$$   a = a_1 \psi_1 + a_2 \psi_2   \bmod N ~~,$$
where $a_1 \in \mathbb Z_p$ and $a_2 \in \mathbb Z_q$ are calculated as
$ ~   a_1   =  a  \hspace{1mm} \bmod p  ~,~  a_2   =  a  \hspace{1mm} \bmod q$.  \\
The four roots $x_1, x_2, x_3, x_4 \in \mathbb Z_N$ of (\ref{main}), represented as positive numbers,
 are obtained using the CRT from the roots $u_1, u_2\in \mathbb Z_p$ and $v_1, v_2\in \mathbb Z_q$
 of the two equations $u^{2} =  C  \hspace{1mm} \bmod p$ and $v^{2} = C  \hspace{1mm} \bmod q$,
 respectively. The roots $u_1$ and $u_2=p-u_1$ are of different parities; likewise, 
 $v_1$ and $v_2=q-v_1$.
If $p$ is congruent $3$ modulo $4$, the root $u_1$ can be computed in deterministic polynomial-time as
$\pm C^{\frac{p+1}{4}} \bmod p$; the same holds for $q$.
If $p$ is congruent $1$ modulo $4$, an equally simple algorithm is not known; however,
 $u_1$ can be computed in probabilistic polynomial-time using Tonelli's algorithm \cite{bach,menezes3} once
 a quadratic non-residue modulo $p$ is known (this computation is the probabilistic part of
 the algorithm), or using the (probabilistic) Cantor-Zassenhaus algorithm \cite{cantor,schipani,jurgen}
 to factor the polynomial $u^2-C$ modulo $p$. 
Using the previous notations, the four roots of (\ref{main}) can be written as  
\begin{equation}
  \label{crt1}
  \left\{ \begin{array}{l} 
       x_1 = u_1 \psi_1 + v_1 \psi_2  \hspace{1cm} \bmod N \\ 
       x_2 = u_1 \psi_1 + v_2 \psi_2  \hspace{1cm} \bmod N \\ 
       x_3 = u_2 \psi_1 + v_1 \psi_2  \hspace{1cm} \bmod N \\ 
       x_4 = u_2 \psi_1 + v_2 \psi_2  \hspace{1cm} \bmod N  ~~. \\ 
     \end{array} \right.
\end{equation}

\begin{lemma}
   \label{lem1}
Let $N=pq$ be a product of two prime numbers. Let $C$ be a quadratic residue modulo $N$; 
 the four roots $x_1, x_2, x_3, x_4$ of the polynomial $x^2-C$ are partitioned
 into two sets $\mathfrak X_1=\{x_1, x_4 \}$ and $\mathfrak X_2=\{x_2, x_3 \}$ such that
 roots in the same set have different parities, i.e. $x_1=1+x_4 \bmod 2$ and $x_2=1+x_3 \bmod 2$.
Furthermore, 
 assuming that $u_1$ and $v_1$ in equation (\ref{crt1}) have the same parity, the residues modulo $p$ and modulo $q$ of each root in $\mathfrak X_1$ have
 the same parity, while each root in $\mathfrak X_2$ has residues of different parities.
\end{lemma}

\noindent
{\sc Proof}.
Since $u_1$ and $v_1$ have the same parity by assumption, then also
$u_2$ and $v_2$ have the same parity. The connection between $x_1$ and $x_4$
 is shown by the following chain of equalities
$$ x_4= u_2 \psi_1 + v_2 \psi_2 =(p-u_1) \psi_1 + (q-v_1) \psi_2 = -x_1 \bmod N = N-x_1 ~~, $$
 because $p\psi_1=0 \bmod N$ and $q\psi_2=0 \bmod N$, and $x_1$ is less than $N$ by assumption, thus
 $-x_1 \bmod N = N-x_1$ is positive and less than $N$. 
A similar chain connects $x_2$ and $x_3=N-x_2$; the conclusion follows because $N$ is odd and thus
$x_1$ and $x_4$ as well as $x_2$ and $x_3$ have different parities.  
\QED
        
\subsection{The Mapping ~~~~$\mathfrak R:~~  x \rightarrow x^2$}
  \label{sect24}
The mapping $\mathfrak R:  ~~ x \rightarrow x^2$ is four-to-one and partitions $\mathbb Z_N^*$
into disjoint subsets $\mathfrak u$ of four elements specified by equation (\ref{crt1}). 
Let $\mathfrak U$ be the group of the four square roots of unity,
that is the roots of $x^2-1$ consisting of the four-tuple
$$  \mathfrak U=\{1, a, -a, -1\}    ~~.           $$
Obviously,  $\mathfrak U$ is a group of order $4$ and exponent $2$.
Each subset $\mathfrak u$, consisting of the four square roots of a given quadratic residue,
 may be described as a coset $m\mathfrak U$ of $\mathfrak U$, i.e.
$$ \mathfrak u= m\mathfrak U=\{m, am, -am, -m\}   ~~.            $$
The number of these cosets is $\frac{\phi(N)}{4}$, and they form a group which is isomorphic
 to a subgroup of $\mathbb Z_N^*$ of order $\phi(N)/4$. 
Once a coset $\mathfrak u=\{ x_1, x_2, x_3, x_4 \}$ is given, 
 a problem is to identify the four elements contained in it. \\
By Lemma \ref{lem1} each $x_i$ is identified by the pair of bits
$$ b_p=(x_i \bmod p) \bmod 2,~~ \mbox{and}~~ b_q=(x_i \bmod q) \bmod 2~~.$$
In summary, the table

\begin{center} 
\begin{tabular}{c|c|c} 
 root   &  $b_p$          &  $b_q$           \\ \hline
 $x_1$  &  $u_1 \bmod 2$  &  $v_1 \bmod 2$   \\
 $x_2$  &  $u_1 \bmod 2$  &  $v_2 \bmod 2$   \\
 $x_3$  &  $u_2 \bmod 2$  &  $v_1 \bmod 2$   \\
 $x_4$  &  $u_2 \bmod 2$  &  $v_2 \bmod 2$   \\
\end{tabular}
\end{center} 

\noindent
shows that two bits identify the four roots. On the other hand, the expression of these two bits involves the prime factorization of $N$, that is $p$ and $q$, but
when the factors of $N$ are not available, it is no longer possible to compute these parity
 bits, and the problem is to find  which parameters can be used, and the minimum number of
 additional bits required to be disclosed in order to label a given root among the four.   \\
Adopting the convention introduced along with equation (\ref{crt1}), a parity bit, namely $b_0\doteq x_i \bmod 2$ distinguishes
 $x_1$ from $x_4$, and $x_2$ from $x_3$, therefore it may be one of the parameters to be used in identifying the four roots.
It remains to determine how to distinguish between roots
having the same parity, without knowing the factors of $N$. 
 
\subsection{Dedekind sums}
  \label{sect25}
A Dedekind sum is denoted by $s(h,k)$ and defined as follows \cite{rademacher}. Let $h,k$ be relatively prime and $k \geq 1$,
 then we set
\begin{equation}
   \label{dede1}
   s(h,k) = \sum_{j=1}^{k} \dede{hj}{k} \dede{j}{k}
\end{equation}   
where the symbol $\deder{x}$, defined as
\begin{equation}
   \label{dede2}
    \deder{x} = \left\{ \begin{array}{ll}
         x-\lfloor x \rfloor - \frac{1}{2}   &  \mbox{if $x$ is not an integer} \\  
         0  & \mbox{if $x$ is an integer}~~,
         \end{array}   \right.
\end{equation}         
denotes the well-known sawtooth function of period $1$. The Dedekind sum satisfies the following
 properties, see \cite{dedekind, grosswald,rademacher} for proofs and details:
\begin{enumerate}
  \item[1)] $h_1 = h_2 \bmod k$  $\Rightarrow$  $s(h_1,k) =s(h_2,k)$
  \item[2)] $s(-h,k) = -s(h,k)$
  \item[3)] $s(h,k)+s(k,h)=- \frac{1}{4}+\frac{1}{12}\left(\frac{h}{k}+\frac{1}{hk}+\frac{k}{h} \right) $, a 
   property known as the reciprocity theorem for Dedekind sums.
  \item[4)] $12 k s(h,k) = k+1- 2 \jacobi{h}{k}  \bmod 8$ for $k$ odd, a property connecting Dedekind sums
             and Jacobi symbols. 
\end{enumerate} 
The first three properties allow us to compute a Dedekind sum by a method that mimics the Euclidean algorithm and has the same efficiency.
In the sequel we need the following Lemma:
\begin{lemma}
   \label{dedelem} 
 If $k=1 \bmod 4$, then, for any $h$ relatively prime with $k$, the denominator of $s(h,k)$ is odd.
\end{lemma}

\noindent
{\sc Proof}. 
In the definition of $s(h,k)$ we can limit the summation to $k-1$ because $\dede{k}{k}=0$,
 furthermore, from the identity $\deder{-x}=-\deder{x}$ it follows that
 $ \sum_{j=1}^{k-1} \dede{hj}{k}=0$ for every integer $h$ \cite{rademacher}, then we may write
$$ s(h,k) = \sum_{j=1}^{k-1} \left(\frac{j}{k}-\frac{1}{2} \right) 
                      \left(\frac{hj}{k} - \left\lfloor \frac{hj}{k} \right\rfloor - \frac{1}{2} \right) =
\sum_{j=1}^{k-1} \frac{j}{k} \left(\frac{hj}{k} - \left\lfloor \frac{hj}{k} \right\rfloor - \frac{1}{2} \right) ~, $$
since $\dede{hj}{k}$ is never $0$, because $j<k$ and $h$ is relatively prime with $k$ by hypothesis.
 The last summation can be split into the sum of two further summations, such that \\
- the first summation $ \sty  \sum_{j=1}^{k-1}\frac{j}{k} \left(\frac{hj}{k} - \left\lfloor \frac{hj}{k} \right\rfloor  \right)  $  has the denominator patently odd; \\
- the second summation is evaluated as $ \sty -\frac{1}{2} \sum_{j=1}^{k-1}\frac{j}{k} = -\frac{k-1}{4} $. 

In conclusion, the denominator of $s(h,k)$ is odd because $s(h,k)$ is the sum of a fraction with odd denominator  with $-\frac{k-1}{4} $, which is an integer number by hypothesis.
\QED

\section{Rabin scheme: primes $p \equiv q \equiv 3 \bmod 4$}
 \label{sect3}

As was said in the introduction, an important issue in using the Rabin scheme is the choice of the right root at the decrypting stage. If $p \equiv q \equiv 3 \bmod 4$, a solution to the identification problem
 has been proposed by Williams \cite{williams} and is reported below, slightly modified from \cite{seberry}, along with three different solutions.
  
\subsection{Williams' scheme}
Williams \cite{seberry,williams}  proposed an implementation of the Rabin cryptosystem, using
 a parity bit and the Jacobi symbol. 

The decryption process is based
 on the observation that, setting 
$D=\frac{1}{2}( \frac{(p-1)(q-1)}{4} +1)$, if $b=a^2 \bmod N$ and $\jacobi{a}{N}=1$, we have
$  b^D =a\jacobi{a}{p}=a \jacobi{a}{q}$, given that 
$$ 
a^{\frac{\varphi(N)}{4}}=(a\psi_1+a\psi_2)^{\frac{\varphi(N)}{4}}=a^{\frac{\varphi(N)}{4}}\psi_1+a^{\frac{\varphi(N)}{4}}\psi_2=\jacobi{a}{p}\psi_1+\jacobi{a}{q}\psi_2=\jacobi{a}{p}=\jacobi{a}{q},  
$$
as
 $a^{\frac{p-1}{2}}=\jacobi{a}{p}\bmod p$, $a^{\frac{q-1}{2}}=\jacobi{a}{q}\bmod q$, and $\frac{p-1}{2}$ and $\frac{q-1}{2}$ are odd (cf. also Lemma 1 in \cite{williams}). 

\begin{description}
  \item[Public-key:]  $[N,S]$, where $S$ is an integer
such that $\jacobi{S}{N}=-1$. 
  \item[Encrypted message]  $   [C,c_1,c_2], $
where 
$$ c_1= \frac{1}{2} \left[ 1- \jacobi{m}{N} \right]  ~~~~,~~~~ \bar m =S^{c_1} m \bmod N   ~~~~,~~~~ 
 c_2=\bar m \bmod 2  ~~~~,\mbox{ and}~~~~   C=\bar m^2 \bmod N    ~~. $$

  \item[Decryption stage]: \\
compute $m'=C^D \bmod N$ and $N-m'$, and choose the number, $m''$ say, with the parity specified by $c_2$.
  The original message is recovered as
$$  m = S^{-c_1} m''  ~~.  $$
\end{description}      
%\vspace{5mm}

\subsection{A second scheme: Variant I}
%We recall here a variant, again exploiting the Jacobi symbol, but in a different way.
%
A simpler variant exploiting the Jacobi symbol is the following: 
\begin{description}
  \item[Public-key:]  $[N]$.
  \item[Encrypted message]  $ [C,b_0,b_1]$, where
$$ C=m^2 \bmod N ~~~~,~~~~  b_0=m \bmod 2  ~~~~\mbox{ and}~~~~
      b_1= \frac{1}{2} \left[ 1+ \jacobi{m}{N} \right] ~~. $$

  \item[Decryption stage]: \\
-  compute, as in (\ref{crt1}), the four roots, written as positive numbers, \\
-  take the two roots having the same parity specified by $b_0$, say $z_1$ and $z_2$, \\
-  compute the numbers
$$  \frac{1}{2} \left[1+\jacobi{z_1}{N}\right] \hspace{10mm} \frac{1}{2} \left[1+\jacobi{z_2}{N}\right]  $$
and take the root corresponding to the number equal to $b_1$. 

\end{description}      

\noindent
{\bf Remark 1.}
The two additional bits are sufficient to uniquely identify $m$ among the four roots,
 because, as previously observed, the roots have the same parity in pairs, and within each of these pairs the roots have opposite Jacobi symbols modulo $N$. In fact, roots with the same parity are of the form $a_1 \psi_1 +a_2 \psi_2$ and $a_1 \psi_1 -a_2 \psi_2$ (or $-a_1 \psi_1 +a_2 \psi_2$), whence the conclusion follows from 
\begin{equation}
\label{jacobi}
 \jacobi{a}{N} = \jacobi{a_1 \psi_1 + a_2 \psi_2 }{pq} =
   \jacobi{a_1 \psi_1 + a_2 \psi_2 }{p }\jacobi{a_1 \psi_1 + a_2 \psi_2 }{ q} =
   \jacobi{a_1 }{p }\jacobi{a_2}{ q} ~~ 
\end{equation}
and
the fact that $-1$ is a nonresidue modulo a Blum prime. 

\subsection{A second scheme: Variant II}
There is a second variant exploiting the Jacobi symbol which, at some extra computational cost
 and further information in the public key,
 requires the delivery of no further bit, since the information needed to decrypt it
 is carried by the encrypted message itself \cite{freeman}.
Let $\xi$ be an integer such that $\jacobi{\xi}{p}=-\jacobi{\xi}{q}=1$, 
 for example $\xi=\alpha^2 \psi_1-\psi_2 \bmod N$, with $\alpha \in \mathbb Z_N^*$. 
The detailed process consists of the  following steps 
\begin{description}
  \item[Public-key:]  $[N,\xi]$.
  \item[Encrypted message]  $ [C]$, where $C$ is obtained as follows
$$ C'=m^2 \bmod N ~~,~~  b_0=m \bmod 2  ~~,~~
      b_1= \frac{1}{2} \left[ 1- \jacobi{m}{N} \right] ~~~\mbox{and}~~~C=C' (-1)^{b_1} \xi^{b_0} \bmod N ~. $$
  \item[Decryption stage]: \\
-  compute $d_0=\frac{1}{2} \left[ 1- \jacobi{C}{q} \right]$, and set $C"= C \xi^{-d_0}$ \\
-  compute $d_1= \frac{1}{2} \left[ 1- \jacobi{C}{N}\right]$, and set $C'= C" (-1)^{d_1}$ \\
-  compute, as in (\ref{crt1}), the four roots of $C'$, written as positive numbers, \\
-  take the root identified by $d_0$ and $d_1$ 

\end{description}      

\noindent
{\bf Remark 2.}
Note that the Jacobi symbol $\jacobi{C}{N}$ discloses the message parity to an eavesdropper.

\subsection{A scheme based on Dedekind sums}
Let $m \in \mathbb Z_N$ be the message to be encrypted, with $N=pq$, $p \equiv q \equiv 3 \bmod 4$. 
The detailed process consists of the  following steps:
\begin{description}
  \item[Public-key:]  $[N]$.
  \item[Encrypted message]  $ [C,b_0,b_1]$, where
$$ C=m^2 \bmod N ~~~~,~~~~  b_0=m \bmod 2  ~~~~,\mbox{ and}~~~~
      b_1= s(m,N)  \bmod 2 ~~, $$
    in which, due to Lemma \ref{dedelem}, the Dedekind sum can be taken modulo $2$
    since the denominator is odd.  
  \item[Decryption stage]: \\
-  compute, as in (\ref{crt1}), the four roots, written as positive numbers, \\
-  take the two roots having the same parity specified by $b_0$, say $z_1$ and $z_2$, \\
-  compute the numbers
$$  s(z_1,N) \bmod 2 \hspace{10mm} s(z_2,N) \bmod 2 ~~,  $$
and take the root corresponding to the number equal to $b_1$. 
\end{description}      

The algorithm works because $s(z_1,N) \bmod 2 \neq s(z_2,N) \bmod 2  $ by the following Lemma.

\begin{lemma}
   \label{dedepari} 
 If $k$ is the product of two Blum primes $p$ and $q$, $(x_1,k)=1$, and $x_2=x_1(\psi_1-\psi_2)$, then 
$$s(x_1,k)+s(x_2,k)=1 \bmod 2 ~~.   $$
\end{lemma}

\noindent
{\sc Proof}.  \\ 
By property 4), which compares the value of the Dedekind sum with the value of the Jacobi symbol,  we have
$$ 12 N s(x_1,N) = N+1- 2 \jacobi{x_1}{N}  \bmod 8 ~~~~\mbox{ and}~~~~  12 N s(x_2,N) = N+1- 2 \jacobi{x_2}{N}  \bmod 8; $$ 
summing the two expressions (member by member) %(term by term)?
and taking into account that $N=1 \bmod 4$ we have
$$ 12N(s(x_1,N)+s(x_2,N)) = 2N+2-2\left[ \jacobi{x_1}{N}+\jacobi{x_2}{N}\right]  \bmod 8 ~, $$
since $12N=4 \bmod 8$, $2N=2 \bmod 8$. Now, we showed above that the sum of the two Jacobi symbols is $0$;
then, applying Lemma \ref{dedelem}, we have
$$ 4(s(x_1,N)+s(x_2,N)) = 4  \bmod 8 ~~ \Rightarrow ~~ s(x_1,N)+s(x_2,N) = 1  \bmod 2 ~~, $$
which concludes the proof.  
\QED
   
\section{Root identification for any pair of primes}
 
If $p$ and $q$ are not both Blum primes, identification of $m$ among the four roots of the
 equation $x^2-C$, where $C=m^2 \bmod N$, can be given by the pair $[b_0,b_1]$ where 
$$   b_0=x_i \bmod 2  ~~~~\mbox{and}~~~~  b_1=(x_i \bmod p) +(x_i \bmod q)  \bmod 2  ~~, $$
 as a consequence of Lemma \ref{lem1}. The bit $b_0$ can be computed at the encryption stage without knowing $p$ nor $q$, while $b_1$ requires, in this definition, $p$ and $q$ to be known, and cannot be directly computed knowing only $N$.  
 
\noindent
In principle, a way to obtain $b_1$ is to publish a pre-computed binary list (or table) that has,
 in position $i$, the bit $b_1$ pertaining to the message $m=i$.
This list does not disclose any useful information on the factorization of $N$ because, even if we know that the residues modulo $p$ and modulo $q$ have the same parity, we do not know which parity, and if these residues  have different parities we do not know which is which. 
Although the list makes the task theoretically feasible, its size is of exponential complexity with respect to
 $N$, and thus practically unrealizable.

\noindent
While searching for different ways of obtaining $b_1$, or some other identifying information, several approaches have been investigated:
\begin{itemize}
 \item to define a polynomial function that assumes the values in the above-mentioned list at the corresponding integer positions; unfortunately this solution is not practical, because this polynomial has a degree roughly equal to $N$, and is not sparse; it is thus more complex than the list. 
  \item to extend the method of the previous section, based on quadratic residues,
 to any pair of primes, by using power residues of higher order and more general reciprocity laws; in particular the quartic reciprocity with Gaussian integers will be involved in providing a neat solution for primes congruent to $5$ modulo $8$.
  \item to exploit group isomorphisms; this could also be of practical interest, although not optimal, in that it relies on the hardness of the Discrete Logarithm problem and it may require more bits than the theoretical lower bound of $2$ to be communicated.
\end{itemize}
      
\subsection{Polynomial function}
  
We may construct an identifying polynomial  as an interpolation polynomial,
 choosing a prime $P$ greater than $N$. Actually the polynomial 
$$  L(x) = \sum_{j=1}^{N-1} \left((j \bmod p)+(j \bmod q) \bmod 2\right)(1-(x-j)^{P-1})   $$ 
assumes the value $1$ in $0< m < N$, if the residues of $m$ modulo $p$ and modulo $q$ have different parities,
 and assumes the value $0$ elsewhere. 
 Unfortunately, as said, the complexity of $L(x)$ is prohibitive and makes this function
 useless in practical terms.

\subsection{Residuosity}
In Section \ref{sect3}, the Jacobi symbol, i.e.  the quadratic residuosity,
 was used to distinguish the roots in the Rabin cryptosystem, when $p=q=3 \bmod 4$. 
For primes congruent $1$ modulo $4$, Legendre symbols cannot distinguish
 numbers of opposite sign, therefore quadratic residuosity is no longer sufficient
 to identify the roots.
Higher power residue symbols could in principle do the desired job, but
unfortunately %  that residues of an order $2^k$ higher than $2$ imply
their use is not straightforward and analogous reciprocity laws or multiplicative properties are not always at hand. 

Actually, higher power residues have been used in some generalizations of the Rabin scheme
 working in residue rings modulo non-prime ideals of algebraic number fields.
For instance, residue rings in Eisenstein or Gauss fields
 were considered in \cite{takagi}, and Rabin-like schemes based on encryption rules
 involving powers of the message higher than $2$ were introduced.
This approach however does not address the problem of separating the roots of a quadratic equation in the
 classic Rabin scheme. 
 
% To fix the situation,  
% 
%It is noted that the Rabin scheme may be generalized in various directions, in particular
% considering residue rings modulo ideals in algebraic number fields. For example, 
% in \cite{takagi}, residue rings modulo ideals in Eisenstein or Gauss fields
% were considered. In particular, an encryption rule involving the fourth power of the
% message in $\mathbb Z[i]$ was proposed together with a companion decryption. 
% , where the encryption rule is modified to a quartic power, and in the language of ideals
% in algebraic number theory, the quartic residuosity is used to separate the roots.
% Although, this approach employs the quartic residuosity in $\mathbb Z[i]$, it does not
% help to solve the problem of separating roots in the classic Rabin scheme,
% as instead do the methods that we describe below.
%  defined in the ring of remainders modulo the product of two primes.
%
%\noindent
Before presenting a neat solution of this root identification problem using quartic reciprocity for primes congruent $5$ modulo $8$, we show below difficulties and attempts concerning a general solution for non Blum primes. \\
Let $2^k$ and $2^h$ be the even exponents of $\mathbb Z_p$ and $\mathbb Z_q$, respectively, that is
  $2^k$ strictly divides $(p-1)$ and $2^h$ strictly divides $(q-1)$, and assume that $k \geq h$. Then the rational power residue symbols $x^{\frac{p-1}{2^k}} \bmod p$
  and $x^{\frac{q-1}{2^h}} \bmod q$ can distinguish, respectively, between $u_1$ and $u_2$ and between $v_1$ and $v_2$. As we would like to use $N$ as a modulo, an idea is to multiply the exponents and consider the function $x^{\frac{\phi(N)}{2^{k+h}}} \bmod N$, which would identify $m$ among the $2^{k+h}$ $2^k$-th roots of unity in $\mathbb Z_N^*$. The idea would be to make these roots publicly available and label them, so that the sender of the message can tell which of them corresponds to the message actually sent. 
  There are two problems: first the exponent $\frac{\phi(N)}{2^{k+h}}$ should also be available, but necessarily in some masked form via multiplication by an odd number in order to hide the factors of $N$; but, most importantly among the public $2^k$-th roots of unity we would find the square roots, and in particular $K\doteq\psi_1-\psi_2$. 
 However, the greatest common divisor of $K+1=2\psi_1$ and $N$ yields $q$, and so $N$ would be factored. 

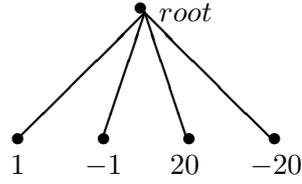
\begin{figure}
\setlength{\unitlength}{0.75mm}
\thicklines
\begin{center}
\begin{picture}(100,40)(0,0)
\put(31.5,33){$\bullet$}
\put(36,32){$root$}
\put(11,11){\line(1,1){22.5}}
\put(26,11){\line(1,3){7.5}}
\put(41,11){\line(-1,3){7.5}}
\put(56,11){\line(-1,1){22.5}}
\put(10,5){$1$}
\put(23,5){$-1$}
\put(38,5){$20$}
\put(52,5){$-20$}
\put(10,10){$\bullet$}
\put(25,10){$\bullet$}
\put(40,10){$\bullet$}
\put(55,10){$\bullet$}
\end{picture}
\end{center}
\caption{Tree representation of the 2-group of order $2\times 2$ in $\mathbb Z_{7\cdot 19}^*$}
\label{fig2}
\end{figure}

Let us now look a bit deeper in this direction, trying to refine this idea.
\noindent
The multiplicative group
 $\mathbb Z_N^*$, direct product of two cyclic groups $\mathfrak C_{p-1}$ and 
 $\mathfrak C_{q-1}$, can also be viewed as the direct product of two abelian subgroups,
 namely a $2$-group and a group of odd order, that is
$$  \mathbb Z_N^* = \left(\mathfrak C_{2^k} \times \mathfrak C_{2^h} \right) \times  
          \left(\mathfrak C_{2f_p+1} \times \mathfrak C_{2f_q+1} \right) ~~.  $$
Therefore, every element $a$ of $\mathbb Z_N^*$ can be written as a product 
 $a_2 a_o$ where $a_o$ is an element of odd order, and $a_2$ is an element of order a power 
 of $2$, i.e. it is an element of a $2$-group which has rank $2$ and exponent $2^k$. \\
The four roots 
$\mathbf V_4=\{1, -1, \psi, -\psi \}$ of $1$, where $\psi = \psi_1 - \psi_2 \bmod N$,
 form a group of order $4$ (the Vierergruppe) of rank $2$, and generators $-1$ and $\psi$. 
Let $a$ be a quadratic residue, then its four square roots $\{A, A_1, A_2, A_3\}$ may be written as $\{A, -A, A\psi, -A\psi\}$, where %$A\psi$
we choose now to consider remainders modulo $N$  of absolute value less than $N/2$. \\
A specific square root $m$ of $a$ among $\{A, -A, A\psi, -A\psi\}$ is identified by
 the sign of $m$ and a further number $c$, possibly a single bit, which should be
 computed with the constraint of using $N$, $m$, and some additional public information
 that should not disclose the factors $p$ and $q$ of $N$. 
Leaving, for a moment, this last constraint, we show how to compute $c$ using a sort of
 residuosity of convenient order depending on the group
  $\mathfrak C_{2^k} \times \mathfrak C_{2^h}$.  \\
Let $2f_N+1=\mbox{lcm}\{2f_p+1,2f_q+1\}$ be the maximum order of the elements in the
 subgroup of odd order, therefore $a_o^{2f_N+1} =1 \bmod N$.
Since $2f_N+1$ and $2^{k}$ are relatively prime, then a generalized Euclidean algorithm
 gives $\alpha$ and $\beta$ such that $\alpha (2f_N+1) +\beta 2^{k}=1$, then we have
$$   m^{\alpha (2f_N+1)} = (m_2 m_o)^{\alpha (2f_N+1)} = m_2^{\alpha (2f_N+1)}
     m_o^{\alpha (2f_N+1)} = m_2^{\alpha (2f_N+1)} = m_2^{1-\beta 2^k} = m_2 ~~. $$
Therefore, an exponentiation with exponent $\alpha (2f_N+1)$ defines an homomorphism
 $\theta$ of the group $\mathbb Z_N^*$ onto the subgroup 
 $\mathcal G_2=\mathfrak C_{2^k} \times \mathfrak C_{2^h}$, 
 such that four-tuples $\mathbf w_b$ of square roots
 of the same element $b$ in $\mathbb Z_N^*$ are mapped into four-tuples 
 $\mathbf g_{\theta(b)}$ of square roots of the same element $\theta(b)$ in $\mathcal G_2$. 
Therefore, in order to identify any specified four-tuple of roots, it is
 sufficient to consider its image in $\mathcal G_2^4$.
Then, it is useful to consider the partition of $\mathcal G_2$ into $4$-tuples that
 are cosets of the group $\mathbf V_4$ of the square roots of $1$. \\
The situation is pictorially described using a $4$-ary rooted tree $\mathfrak T$
 with nodes labeled by the elements of the $2$-group $\mathcal G_2$. 
The  four nodes at the first layer below the root are labeled by the four roots of unity. 
At this layer, the node labeled with $1$ is a terminal node, the remaining
 three nodes may or may not be terminal nodes depending on the form of the primes 
 $p$ and $q$. 
The height of the tree is $k (\geq h)$; the number of nodes at each level is a 
 multiple of $4$, and depends on the forms of the primes $p$ and $q$. 
If there is a path (a sequence of branches) connecting a node $u$ with a node $v$ of
 a superior layer, we say that $v$ is above $u$. \\
For the sake of example, Figures  \ref{fig2},  \ref{fig3}, and \ref{fig1} show every 
 possible shape of trees with at most two layers. In particular, the tree in 
 Figure \ref{fig2} corresponds to a pair of primes congruent $3$ modulo $4$, the tree in 
 Figure \ref{fig3} corresponds to a pair of primes, one congruent $3$ modulo $4$ and
 the second congruent $5$ modulo $8$, lastly, the tree  in Figure \ref{fig1}  
 corresponds to a pair of primes congruent $5$ modulo $8$. \\
Note that every set of four nodes, directly connected to the same node, can be identified
 by a single label, say the coset leader, since the set of labels of these nodes can be seen as a coset
 of $\mathbf V_4$. 

%\vspace{10mm}

\begin{figure}
\setlength{\unitlength}{0.75mm}
\thicklines
\begin{center}
\begin{picture}(185,75)(0,0)
\put(67,70){$root$}
\put(64,70){$\bullet$}
\put(6,41){\line(2,1){60}}
\put(5,40){$\bullet$}
\put(5,43){$1$}
\put(26.5,33.0){\line(1,1){37.5}}
\put(92.5,34){\line(-4,5){28.0}}
\put(156.5,40){\line(-3,1){92.0}}
\put(26.5,33){$\bullet$}
\put(31,32){$-1$}
\put(91.5,33){$\bullet$}
\put(96,32){$6$}
\put(156.5,39){$\bullet$}
\put(160,38){$-6$}
\put(136,17){\line(1,1){22.5}}
\put(151,17){\line(1,3){7.5}}
\put(166,17){\line(-1,3){7.5}}
\put(181,17){\line(-1,1){22.5}}
\put(135,11){$22$}
\put(150,11){$-22$}
\put(165,11){$27$}
\put(180,11){$-27$}
\put(135,16){$\bullet$}
\put(150,16){$\bullet$}
\put(165,16){$\bullet$}
\put(180,16){$\bullet$}
\end{picture}
\end{center}
\caption{Tree representation of the 2-group of order $2^{2}\times 2$ in $\mathbb Z_{5\cdot 7}^*$}
\label{fig3}
\end{figure}
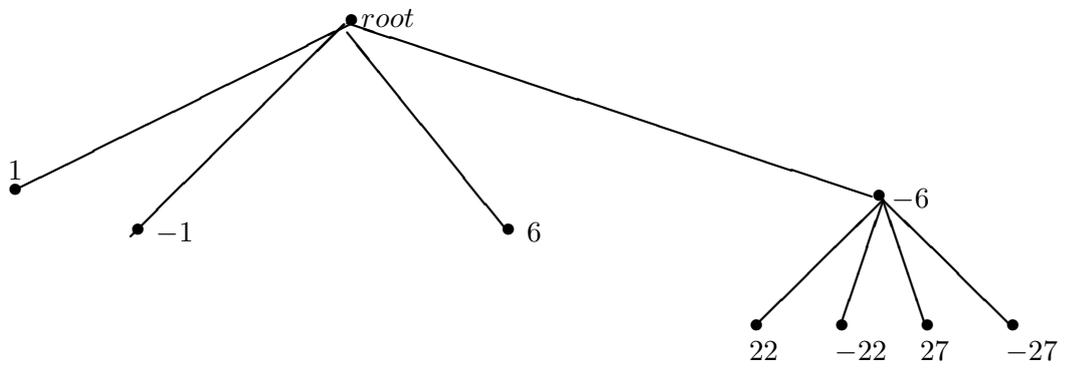

\begin{figure}
\setlength{\unitlength}{0.75mm}
\thicklines
\begin{center}
\begin{picture}(185,75)(0,0)
\put(67,70){$root$}
\put(64,70){$\bullet$}
\put(6,41){\line(2,1){60}}
\put(5,40){$\bullet$}
\put(5,43){$1$}
\put(26.5,33.0){\line(1,1){37.5}}
\put(92.5,34){\line(-4,5){28.0}}
\put(156.5,40){\line(-3,1){92.0}}
\put(26.5,33){$\bullet$}
\put(31,32){$-1$}
\put(6,11){\line(1,1){22.5}}
\put(21,11){\line(1,3){7.5}}
\put(36,11){\line(-1,3){7.5}}
\put(51,11){\line(-1,1){22.5}}
\put(5,5){$8$}
\put(18,5){$-8$}
\put(34,5){$18$}
\put(47,5){$-18$}
\put(5,10){$\bullet$}
\put(20,10){$\bullet$}
\put(35,10){$\bullet$}
\put(50,10){$\bullet$}
\put(91.5,33){$\bullet$}
\put(96,32){$14$}
\put(71,11){\line(1,1){22.5}}
\put(86,11){\line(1,3){7.5}}
\put(101,11){\line(-1,3){7.5}}
\put(116,11){\line(-1,1){22.5}}
\put(70,5){$12$}
\put(83,5){$-12$}
\put(98,5){$27$}
\put(112,5){$-27$}
\put(70,10){$\bullet$}
\put(85,10){$\bullet$}
\put(100,10){$\bullet$}
\put(115,10){$\bullet$}
\put(156.5,39){$\bullet$}
\put(160,38){$-14$}
\put(136,17){\line(1,1){22.5}}
\put(151,17){\line(1,3){7.5}}
\put(166,17){\line(-1,3){7.5}}
\put(181,17){\line(-1,1){22.5}}
\put(135,11){$21$}
\put(150,11){$-21$}
\put(165,11){$31$}
\put(180,11){$-31$}
\put(135,16){$\bullet$}
\put(150,16){$\bullet$}
\put(165,16){$\bullet$}
\put(180,16){$\bullet$}
\end{picture}
\end{center}
\caption{Tree representation of the 2-group of order $2^{2}\times 2^2$ in $\mathbb Z_{5\cdot 13}^*$}
\label{fig1}
\end{figure}
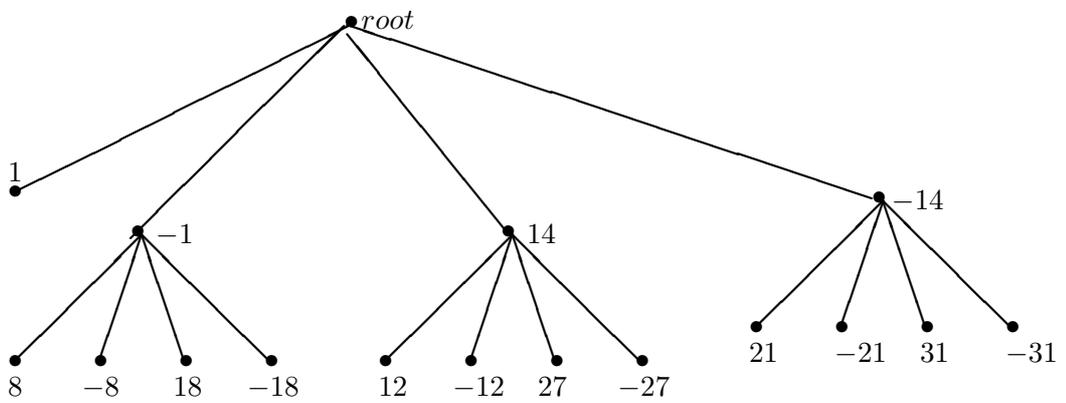

\noindent
The next two lemmas show how to use the tree to identify the correct root of 
 $X^2=b \bmod N$, but, unfortunately, also show that the residuosity connected with
 $\theta(.)$ discloses the factorization of $N$. 
 
\begin{lemma}
   \label{lemsplit}
Assume that the exponent $\alpha (2f_N+1)$ is public
  together with a table $T$ of $2^{h+k-2}$ elements, containing one of the two positive elements for
    each set of $4$ elements of the group $\mathcal G_2$ as described in its tree
    representation.
Then, only two bits are sufficient to identify a square root $m$ of $b$, that is,
 one bit for the sign of $m$, and one bit telling whether $|\theta(m)|$, the absolute value of $\theta(m)$, can be found
 in the table or not.   
\end{lemma}

%\begin{proof}
\noindent
{\sc Proof}
When the sender wants to encrypt $m$, then the triple $\{b,b_0,b_1\}$ is sent, where  $b=m^2 \bmod N$,
 $b_0$ is the sign of $m$, and $b_1=\mathfrak{I}(|\theta(m)|\in T)$, with $\mathfrak{I}$ being the indicator function. \\
Given $[b,b_0,b_1]$ and knowing the factorization of $N=pq$, the right value $m$ is
 identified as follows: 
\begin{enumerate}
  \item Solve the equation $x^2=b \bmod N$ and find four values $[A, -A, B, -B]$
  \item Compute $[|\theta(A)|,|\theta(B)|]$, one of these two values is in the table,
    therefore select the one compatible with $b_1$.
  \item Define the correct value $m$ using the previous value and $b_0$.
\end{enumerate}
\QED %\end{proof}

\noindent
Unfortunately, the disclosure of $\alpha (2f_N+1)$ leads to factor $N$. 

\begin{lemma}
   \label{lemprob}
Assuming that $\alpha (2f_N+1)$ is known, then the probability of factoring $N$
 is not less than $1/2$.
\end{lemma}

%\begin{proof}
\noindent
{\sc Proof}
We already showed that knowing $\psi$ we can factor $N$.
Picking an integer $x_r$ at random, the probability that $u=x_r^{\alpha (2f_N+1)}$ is below $\psi$ or $-\psi$ in the tree is at least $1/2$. In the favorable event that $u$ is below $\psi$ or $-\psi$, a power of $u$
 with a convenient exponent $2^{f(u)}$ gives $\psi$. The probability is exactly $1/2$
 in the case of Blum primes, otherwise is larger as can be deduced from the trees.
\QED %\end{proof}

\noindent
In conclusion, the scheme allows us in principle to compute two bits discriminating the four
 roots of $b$, by means of functions computable using
 only $m$, $N$ and not its
 factorization. Unfortunately, the additional information made public, the table and the exponent, permit the 
 factorization of $N$ deterministically, as one can retrieve $\psi$ from the table,  as well 
 as probabilistically with high probability, as a consequence of Lemma \ref{lemprob}. \\
Therefore, it is necessary to look at different kinds of higher order residuosity which
 should allow
\begin{itemize}
  \item a definition of symbols a la Jacobi specifying the residue character;
  \item a reciprocity law for these symbols;
  \item the values of the symbols should belong to a finite group which does not reveal any
   information allowing the factorization of $N$.
\end{itemize}   

\noindent
Let $\ell$ denote the height of the tree $\mathfrak T$, and $\zeta_{2^{\ell}}$ be
 a primitive root of unity;
 it turns out that such a $2^{\ell}$-residuosity exists in the ring of integers 
 $\mathbb Z[\zeta_{2^{\ell}}]$ of cyclotomic fields $\mathbb Q(\zeta_{2^{\ell}})$.
Let $\nu \in \mathbb Z[\zeta_{2^{\ell}}]$ be irreducible.
A symbol of residuosity may be defined, \cite[Theorem 46, p.211]{taylor}, as
\begin{equation}
  \label{gaussjacsymb}
 \gauss{b}{\nu}_{2^{\ell}} = b^{\frac{\mathcal N(\nu)-1}{2^{\ell}}} \bmod \nu =
 \zeta_{2^{\ell}}^{\gamma(b)} , %~~\forall b \in \mathbb Z[\zeta_{2^{\ell}}] ~~ ,
\end{equation}
where $\mathcal N(\nu)$ is the norm of $\nu$ in $\mathbb Q(\zeta_{2^{\ell}})$,
 $\gamma(b)$ is an integer that certainly exists, since $\zeta_{2^{\ell}}$ and $b^{\frac{\mathcal N(\nu)-1}{2^{\ell}}}$
 are both roots of $X^{2^{\ell}}-1 \bmod \nu$. \\
Using this residuosity, the difficulty is moved to compute $\gamma(b)$, however,
 in case of quartic residuosity, the task is made possible by the
 Gauss-Jacobi's quartic residue symbols and their reciprocity law, as we show in the next subsection. 

\subsubsection{Identification scheme using quartic residuosity}

Assuming that $p$ and $q$ are congruent $5$ modulo $8$,  we show here that the quartic residuosity in the Gaussian integers
 is sufficient to discriminate the $4$ square roots with exactly $2$ bits. 

Let $\mathbb Z[i]$ be the ring of Gaussian integers, which is Euclidean, so that the
 factorization is unique except for a reordering and a multiplication by units.
The units are $\mathfrak U=\{1, -1, i, -i\}$ and form a cyclic group \cite{hardy}. 
Any integer $z=x+iy$ in $\mathbb Z[i]$ has four associates, namely $z, -z, iz$, and $-iz$.
In $\mathbb Z[i]$ the rational primes $p$ congruent $1$ modulo $4$ split as 
 $p=(a+ib)(a-ib)$, and $2$ splits as $2=(1+i)(1-i)$. 
The following notions and properties are taken from 
 \cite[p.119-127]{rosen}, which we refer to for proofs and details.
 
\begin{definition}
An integer $x+iy$ of $\mathbb Z[i]$, is said to be {\em primary} if $x+iy=1 \bmod (1+i)^3$. \\
An integer $z \in \mathbb Z[i]$ is said to be odd if it is not divisible by $1+i$. \\
The norm of $x+iy \in \mathbb Z[i]$ is $\mathcal N(x+iy)= x^2+y^2$.
\end{definition}

\noindent
We note that any odd integer $x+iy$ has an associated primary which can be obtained
  upon multiplication by a unit. 
We now prove the following theorem:
 
\begin{theorem}
  \label{pmod8}
An odd prime $p$ congruent $5$ modulo $8$ has the representation, as a sum of two squares,
 of the form  $p=(2X+1)^2+4(2Y+1)^2$, then in $\mathbb Z[i]$ decomposes as
 $$ p =((2X+1)+2(2Y+1)i) ((2X+1)-2(2Y+1)i) ~~, $$
and a primary factor is
$$ \pi =  ((2X+1)+2(2Y+1)i) (-1)^{X-1} ~~. $$  
\end{theorem}

%\begin{proof}
\noindent
{\sc Proof}
Since $p$ can be written as a sum of two squares $p=(2X+1)^2+4y^2$, the first part of the lemma is proved by showing that $y$ is odd.
Taking $p$ modulo $8$ we have
$$ p \bmod 8=5 = 4X(X+1)+1+4y^2 = 1+4y^2 \Rightarrow 4 y^2 =4 \bmod 8  \Rightarrow y^2 =1 \bmod 2 ~~, $$
which implies $y=1 \bmod 2$. \\
The prime factor $\pi$ of $p$ in $\mathbb Z[i]$ is primary if it is congruent $1$ modulo
 $-2+2i$. Imposing this condition, with $u$ a unit, and considering that $4=0 \bmod (-2+2i)$,
 we have
$$ 1 = ((2X+1)+2i(2Y+1)) u \bmod (-2+2i) = u (2X+3) \bmod (-2+2i)~~,  $$ 
because $2i=2  \bmod (-2+2i)$. We distinguish two cases:
\begin{enumerate}
  \item If $X$ is even then $u$ must satisfy the condition $3u=1 \bmod (-2+2i)$, which forces
     $u=-1$, that is $u=(-1)^{X-1}$.
  \item If $X$ is odd then $u$ must satisfy the condition $5u=1 \bmod (-2+2i)$, which forces
     $u=1$, that is $u=(-1)^{X-1}$ again.
\end{enumerate}
This concludes the proof.     
\QED %\end{proof}

\noindent
Let $\pi \in \mathbb Z[i]$ be an odd irreducible, and $\pi \notdiv \alpha$. There exists a
 unique integer $j$, \cite[p.122]{rosen}, such that
$$  \alpha^{\frac{\mathcal N(\pi)-1}{4}} = i^j  \bmod \pi ~~.  $$
This property is used to define a quartic residue symbol as
$$ \gaussf{\alpha}{\pi} = \left\{ \begin{array}{ll}
         i^j & \mbox{if} ~~\pi \notdiv \alpha \\
          0  & \mbox{otherwise} 
     \end{array}  \right.  ~~.  $$     
Let $\nu = a+ib$ be a primary odd number, then a Jacobi-like symbol for quartic residues, 
 written as $\gaussf{\beta}{\nu}$ and called Gauss-Jacobi symbol, is defined multiplicatively, similarly to 
 the Jacobi symbol in the quadratic case.
It satisfies the following properties \cite{rosen, lemmer} that 
 allow us to evaluate the symbol without knowing the factorization of the arguments:
\begin{enumerate}
   \item  $\gaussf{\alpha+\mu \nu}{\nu}= \gaussf{\alpha}{\nu}$, %see \cite{rosen}
   \item  $\gaussf{\alpha\beta}{\nu}= \gaussf{\alpha}{\nu}\gaussf{\beta}{\nu}$, %see \cite{rosen}
   \item  
   $\gaussf{i}{\nu}= i^{-\frac{a-1}{2}}$ and thus 
          $\gaussf{-1}{\nu}= (-1)^{\frac{a-1}{2}}$, %see \cite{lemmer}
   \item  $\gaussf{1+i}{\nu}= i^{\frac{a-1-b-b^2}{4}}$ and thus
          $\gaussf{2}{\nu}= i^{\frac{-b}{2}}$, %see \cite{lemmer}
   \item  If $\omega = c + di$ is a primary odd number, its real part $c$ is odd,
     then either $c$ or $-c$ is congruent $1$ modulo $4$, it follows that the real part
     of $\omega$ or $-\omega$ is congruent $1$ modulo $4$.
     Let $\alpha = u + vi$ and $\beta=t+wi$ be odd with the real part congruent $1$
     modulo $4$, the reciprocity law takes the Jacobi-Kaplan form  %\cite{rosen,lemmer}
$$  \gaussf{\alpha}{\beta} \gaussf{\beta}{\alpha}^{-1} = (-1)^{\frac{v\cdot w}{4}}   ~~.   $$
\end{enumerate}

\noindent
The main theorem of this section permits to identify the four square roots of 
 a quadratic residue in $\mathbb Z[i]$ using only two bits and without unveiling the factorization of $N$ or $\nu$. \\
Proceeding as we previously did in the definition of $\psi_1$ and $\psi_2$, given $\pi_1, \pi_2$ relatively prime integers in $\mathbb Z[i]$, we can find $\xi_1,\xi_2$, such that 
 $\xi_1 +\xi_2=1$,  $\xi_1=\mu_2 \pi_2$ and  $\xi_2=\mu_1 \pi_1$ for some $\mu_1, \mu_2\in \mathbb Z[i]$, and set $\xi=\xi_1 -\xi_2$, which turns out to be a square root of $1$ modulo $\nu$, i.e. $\xi^2=1 \bmod \nu$.

\begin{theorem}
  \label{maingauss}
Let $\nu = \pi_1 \pi_2$ be the product of two primary primes having norms congruent $5$ modulo $8$. A root $\alpha$ among the four square roots 
$\{\gamma, -\gamma, \gamma\xi, -\gamma\xi \}$ of a quadratic residue $\beta$ can be uniquely identified with two bits $b_0$ and $b_1$ defined as:
$$  b_0 = \left\{ \begin{array}{lcl}
                  1  & \mbox{if} & \Re(\alpha) > 0 \\     
                  0  & \mbox{if} & \Re(\alpha) < 0 \\     
                 \end{array}  \right.  ~~, $$
(use  $\Im(\alpha)$ if $\Re(\alpha)=0$); 
$$  b_1 =  \left\{ \begin{array}{lcl}
                  1  & \mbox{if} & \gaussf{\alpha}{\nu} \in \{1,i\} \\     
                  0  & \mbox{if} & \gaussf{\alpha}{\nu} \in \{-1,-i\} \\     
                 \end{array}  \right.   %\frac{1+\gaussf{\alpha}{\nu}}{2} 
~~. $$
\end{theorem}

%\begin{proof}
\noindent
{\sc Proof}.
With the given choice of $b_0$, the parameter $b_1$ must discriminate $\alpha$ from
 $\alpha\xi$ or $-\alpha\xi$. Due to the multiplicative property of the Gauss-Jacobi
 symbol, this is tantamount to show that $\gaussf{\xi}{\nu}=\gaussf{-\xi}{\nu}=-1$. 
 Indeed we have
$$ \gaussf{\xi}{\nu}= \gaussf{\xi}{\pi_1} \gaussf{\xi}{\pi_2}=  \gaussf{\xi_1 -\xi_2}{\pi_1} \gaussf{\xi_1 -\xi_2}{\pi_2} = \gaussf{\xi_1}{\pi_1} \gaussf{ -\xi_2}{\pi_2}  ~~. $$
But $\xi_1=1-\xi_2$ and conversely $\xi_2=1-\xi_1$, so we obtain the expression
$$ \gaussf{\xi}{\nu}= \gaussf{1-\xi_2}{\pi_1} \gaussf{ -(1-\xi_1)}{\pi_2} =
   \gaussf{1}{\pi_1} \gaussf{ -1}{\pi_2} = -1 ~~. $$
This conclusion follows because, by Theorem \ref{pmod8}, $\pi_2$ is of the form $((2X+1)+2(2Y+1)i) (-1)^{X-1} $, which implies
$$ \gaussf{ -1}{\pi_2} = (-1)^{\frac{(2X+1)(-1)^{X-1}-1}{2}}= -1 ~~,  $$
since the exponent is always odd, whatever be the parity of $X$. \\
In the same way $\gaussf{-\xi}{\nu}=-1$ by exchanging the role of $\pi_1$ and $\pi_2$. \\
In summary, $\gaussf{\alpha}{\nu}=-\gaussf{\alpha\xi}{\nu}=-\gaussf{-\alpha\xi}{\nu}$, then $b_1\in\{0,1\}$ distinguishes among the two roots  with the same $b_0$. 
\QED %\end{proof}

\paragraph{Remark 3.}
Let $N$ be equal to the norm of $\nu$, then a representation of the elements
 of the finite ring $\mathfrak Z_{\nu}=\mathbb Z[i]/\nu\mathbb Z[i]$, which is 
 isomorphic to $\mathbb Z_N$, may consist of the same elements of $\mathbb Z_N$. 
A more "natural" representation of $\mathfrak Z_{\nu}$ consists of $N$ elements
 of $\mathbb Z[i]$, which have minimum Euclidean norm and are not congruent to
 one another modulo $\nu$. 
The two representations are perfectly equivalent, the use of
 one or the other only depends on the simplicity of computations and arithmetic 
 operations in $\mathbb Z[i]/\nu\mathbb Z[i]$.

\vspace{3mm}
\noindent
Let $N=pq$ be decomposed in $\mathbb Z[i]$ as a product $\nu \bar \nu$, where
 $\nu= \pi_1 \pi_2$ is the product of an irreducible factor of $p$ and an irreducible 
 factor of $q$.
Noting that $\nu$ and $\bar \nu$ are relatively prime, any number $f$ of $\mathbb Z_N$
 is uniquely identified by the pair $[f_1,f_2]$ obtained by taking the remainders modulo
 $\nu$ and modulo $\bar \nu$, i.e. $f_1= f \bmod \nu$, $f_2= f \bmod \bar \nu$, and $f_2$ 
 is easily seen to be the complex conjugate $\bar{f_1}$ of $f_1$. 
The value $f$ is recovered from the pair $[f_1,\bar{f_1}]$, by using the Chinese remainder
 theorem 
\begin{equation}
  \label{crtcpx}
  f = f_1 \zeta_1  + \bar{f_1} \zeta_2 ~~\bmod N ~~,
\end{equation}
 where $\zeta_1=\mu_1 \bar \nu \bmod N$ and  $\zeta_2=\mu_2 \nu \bmod N$, are the complex
 counterpart of $\psi_1$ and $\psi_2$, with $\mu_1$ and $\mu_2$ 
 computed by means of the generalized Euclidean algorithm.
% , which is started  to compute the $\gcd$ between $\nu$ and $\bar \nu$.
It is pointed out, as a consequence of equation (\ref{crtcpx}), that a quadratic residue
 $m$ modulo $N$ is also a quadratic residue modulo $\nu$, and a square root $A$ of
 $m$ modulo $N$ corresponds to a square root $\alpha$ of $m$ modulo $\nu$. 
Therefore, $-A$ corresponds to $-\alpha$, $A\psi$ corresponds to
  $\alpha\xi$, and  $-A\psi$ corresponds to $-\alpha\xi$ because $\xi = \psi \bmod \nu$.
This last identity is straightforwardly proved observing that 
$$ 1=\psi_1+\psi_2= \lambda_2 q +\lambda_1 p= \lambda_2 \pi_2 \bar \pi_2 +\lambda_1 \pi_1 \bar \pi_1 ~~ $$
in $\mathbb Z[i]$, thus, taking the remainder modulo $\nu$, we have
$$ 1= (\lambda_2 \bar \pi_2) \pi_2 +(\lambda_1 \bar \pi_1) \pi_1 \bmod \nu =
      \xi_1 + \xi_2 \bmod \nu    ~~, $$
due to the definition of $\xi_1$ and $\xi_2$, and finally
 $\xi_1=\psi_1 \bmod \nu$ and $\xi_2=\psi_2 \bmod \nu$ for the Chinese remainder theorem.
   
\vspace{3mm}
\noindent 
A Rabin scheme working with primes $p$ and $q$ congruent $5$ modulo $8$ can be defined
 considering the decomposition $N=\nu \bar \nu$ with $\nu=\pi_1 \pi_2$
 being the product of two primary factors of $p$ and $q$ respectively.
 
\begin{description}
  \item[Public-key:]  $[\nu]$.
  \item[Message:]  $m$.
  \item[Encrypted message]  $ [C,b_0,b_1]$, where
$$ C=m^2 \bmod N ~~~~,~~~~  b_0=m \bmod 2  ~~~~,\mbox{ and}~~~~
      b_1 =  \left\{ \begin{array}{lcl}
                  1  & \mbox{if} & \gaussf{m}{\nu} \in \{1,i\} \\     
                  0  & \mbox{if} & \gaussf{m}{\nu} \in \{-1,-i\} \\     
                 \end{array}  \right.  ~~. $$ 
  \item[Decryption stage]: \\
-  compute, as in (\ref{crt1}), the four roots of $C$ modulo $N$, written as positive numbers, \\
-  take the two roots having the same parity specified by $b_0$, say $z_1$ and $z_2$, \\
-  compute the quartic residues 
$$ \gaussf{z_1}{\nu} \hspace{10mm} \gaussf{z_2}{\nu}  ~~,  $$
and take the root corresponding to $b_1$. 
\end{description}

\paragraph{Remark 4.}
The extension to residuosity of higher order is straightforward only up to 
 $\mathbb Q(\zeta_{32})$
 because these fields are Euclidean \cite{kaiblinger}. The next field 
 $\mathbb Q(\zeta_{64})$ has class number $17$, thus is certainly not Euclidean.
Also, the Euclide algorithm may not always be easy to perform.
It is known that for Gaussian integers $\mathbb Z[i]$
 the division may be performed by {\em rounding} the entries of the quotient of
 integers $v = \frac{v_0+v_1i}{\nu} = (a_0; a_1) \in \mathbb Q^2$ to nearest integers,
 $a'_0 +a'_1 i= ( \lfloor a_0 + \frac{1}{2} \rfloor; \lfloor a_1 + \frac{1}{2} \rfloor) \in \mathbb Z^2$. The remainder of minimum norm is obtained as 
 $r_0+r_1i =(v_0+v_1i)-\nu \cdot (a_0 +a_1 i)$.      \\
It is also known that the 1-step norm-Euclidean algorithm for $\mathbb Z[\zeta_8]$ devised by
 Eisenstein \cite{eisenstein} is implicitly defined by {\em rounding}, and 
 \cite[Sec. 4.1]{monico} includes an explicit proof.

\subsection{Group isomorphisms}
In this section, we describe a practical method, working with any pair of primes, that can
 have acceptable complexity, although it requires a one-way function that might be weaker
 than factoring.
 
A possible solution is to use  a function $\mathfrak d$ defined from $\mathbb Z_N$ into a group $\mathfrak G$ of the same order, 
 and define a function $\mathfrak d_1$ such that $\mathfrak d_1(x_1)=\mathfrak d(x_2)$.
 The public key consists of the two functions $\mathfrak d$ and $\mathfrak d_1$. At the encryption stage,
 both are evaluated at the same argument, the message $m$, and the minimum information necessary to
 distinguish their values is delivered
 together with the encrypted message. The decryption operations are obvious.
The true limitation  of this scheme is that $\mathfrak d$ must be a one-way function, 
otherwise two square roots that allow us to factor $N$ can be recovered as in the residuosity subsection.

Following this approach, we propose the following solution, based on the hardness of computing discrete logarithms. 

Given $N$, let $P=\mu N +1$ be a prime (the smallest prime), that certainly exists by Dirichlet's theorem
 \cite{apostol}, that is congruent $1$ modulo $N$. Let $g$ be a primitive element generating the multiplicative group $\mathbb Z_P^*$.
 
Define $g_1=g^{\mu}$ and $g_2=g^{\mu(\psi_1-\psi_2)}$, and as usual let $m$ denote the message. 

\begin{description}
\item[Public key:] $[N,P,g_1,g_2]$. 
\item[Encryption stage:]  $[C, b_0, d_1, d_2, p_1, p_2]$, where $C=m^2 \bmod N$, $b_0=m \bmod 2$,
  $p_1$ is a position in the binary expansion of $g_1^{m} \bmod P$, whose bit $d_1$ is different
  from the bit in the  corresponding position of the binary expansion of $g_2^{m} \bmod P$, 
  and $p_2$ is a position in the binary expansion of $g_1^{m} \bmod P$, whose bit $d_2$ is different
  from the bit in the corresponding position of the binary expansion of $g_2^{-m} \bmod P$. \item[Decryption stage]: \\
-  compute, as in (\ref{crt1}), the four roots, written as positive numbers, \\
-  take the two roots having the same parity specified by $b_0$, say $z_1$ and $z_2$, \\
-  compute $A=g_1^{z_1} \bmod P$ and $B=g_1^{z_2} \bmod P$ \\
-  between $z_1$ and $z_2$, the root is selected that has the correct bits $d_1$ and $d_2$ in both
  the given positions $p_1$ and $p_2$ of the binary expansion of $A$ or $B$. 
\end{description}

The algorithm is justified by the following Lemma.

\begin{lemma}
  \label{lemmaprag}
The power $g_0=g^{\mu}$ generates a group of order $N$ in $\mathbb Z_P^*$, thus the correspondence
 $x \leftrightarrow g_0^x $ establishes an isomorphism between a multiplicative subgroup of $\mathbb Z_P^*$
 and the additive group of $\mathbb Z_N$. The four roots of $x^2=C \bmod N$, $C=m^2 \bmod N$ are in a one-to-one correspondence
 with the four powers $g_0^m \bmod P$, $g_0^{-m} \bmod P$, $g_0^{m(\psi_1-\psi_2)} \bmod P$ and 
 $g_0^{-m(\psi_1-\psi_2)} \bmod P$.
\end{lemma}

\noindent
{\sc Proof}.
The first part is due to the choice of $P$: the group generated by $g_0$ has order $N$, thus,
 the isomorphism follows immediately.
The second part is a consequence of Section \ref{sect24}.
\QED

The price to pay is the costly arithmetic in  $\mathbb Z_P$, and the equivalence of the security of the
 Rabin cryptosystem
 with the hardness of factoring is now conditioned by the complexity of computing the discrete logarithm in $\mathbb Z_P$.

\section{The Rabin signature}
   \label{sect21}

In the introduction, we said that a Rabin signature of a message $m$ may consist of a pair
 $[n,S]$; however, if $x^2 =m \bmod N$ has no solution, this signature cannot be directly
 generated. To overcome this obstruction, a random pad $U$ was proposed \cite{seberry}, and attempts
 are repeated until $x^2 =mU \bmod N$ is solvable, and the signature is the triple
 $(m,U,S)$, \cite{seberry}. A verifier compares $mU \bmod N$ with $S^2$ and accepts the signature as valid
 when these two  numbers are equal. 

This section presents a modified version of this scheme, where $U$ is computed deterministically. 

Now, the quadratic equation $x^2=m \bmod N$ is solvable if and only if $m$ is a quadratic residue
 modulo $N$, that is $m$ is a quadratic residue modulo $p$ and modulo $q$.
When $m$ is not a quadratic residue, 
 we show below how to exploit the Jacobi symbol to compute a suitable pad and obtain quadratic residues modulo $p$ and $q$.
Let \\
$f_1 = \sty \frac{m_1}{2}\left[1-\jacobi{m_1}{p}\right]+ \frac{1}{2}\left[1+\jacobi{m_1}{p}\right]$ ~,  ~ $f_2 = \sty \frac{m_2}{2} \left[1-\jacobi{m_2}{q}\right]+ \frac{1}{2} \left[1+\jacobi{m_2}{q}\right]$. %, with $n_1$ and $n_2$ two quadratci nonresidues with respect to $p$ and $q$, respectively.
 
\noindent
Writing $m=m_1 \psi_1 +m_2 \psi_2$, the equation
$$   x^2 =(m_1 \psi_1 +m_2 \psi_2) (f_1 \psi_1 +f_2\psi_2) =
     m_1f_1 \psi_1 +m_2 f_2\psi_2$$
is always solvable modulo $N$, because $m_1f_1$ and $m_2 f_2$
 are clearly quadratic residues modulo $p$ and modulo $q$, respectively, since $\jacobi{m_1}{p}=\jacobi{f_1}{p}$, $\jacobi{m_2}{q}=\jacobi{f_2}{q}$, so that
$$   \jacobi{m_1 f_1}{p} = \jacobi{m_1}{p}\jacobi{f_1}{p} = 1 ~,  ~ \jacobi{m_2 f_2}{q} = \jacobi{m_2}{q}\jacobi{f_2}{q} = 1
  ~~. $$

Note that if $p$ and $q$ are Blum primes, it is possible to choose $f_1=  \jacobi{m_1}{p}$ and $f_2=  \jacobi{m_2}{q}$. 

Thus we can describe the following procedure:
\begin{description}
  \item[Public-key:] $N$ 
  \item[Signed message:] $   [U,m,S], $ where 
    $U= R^2\left[f_1\psi_1+f_2\psi_2\right] \bmod N$ is the padding factor,
    with $R$ a random number, and $S$ is any solution of the equation $x^2=mU \bmod N$. $R$ is needed to avoid that knowing $U$ allows to easily factor $N$.
  \item[Verification:] 
 compute $mU \bmod N$ and $S^2 \bmod N$; the signature is valid if and only if these two numbers are equal.
\end{description}  

This signature scheme has several interesting features: 

\begin{enumerate}
\item the signature is possible using every pair of primes, and thus it could be used with 
   the modulo of any RSA public key, for example;
\item different signatures of the same document are different;
\item the verification needs only two multiplications, therefore it is fast enough to be used  in authentication protocols. 
\end{enumerate}  

\subsection{Forgery attacks}

Schemes of this type are however vulnerable to forgery attacks: it is relatively easy to compute $S^2 \bmod N$, choose any message $m'$, compute $U'=S^2 m'^{-1} \bmod N$, and forge the signature as $(m',U',s)$ without knowing the factorization of $N$. In some variants a hash $H(m)$ is used instead of $m$ and $S$ is a solution of $x^2 = H(mU) \bmod N$, but this does not help against the above forgery attack.
The following variant aims at countering this vulnerability.

\begin{description}
  \item[Public-key:] $N$ 
  \item[Signed message:] $   [m,UK^2\bmod N,SK^3\bmod N, K^4\bmod N], $ where 
    $U$ is the padding factor,
    $K$ a random number, and $S$ is any solution of the equation $x^2=mU \bmod N$.
  \item[Verification:] 
 compute $(SK^3)^2 \bmod N$ and $mUK^2K^4 \bmod N$; the signature is valid if and only if these two numbers are equal.

\end{description}  

We remark that $U$, $K$ and $S$ are not known. Forgery would be possible if $K$ were known, but to know $K$ one has to solve an equation of degree at least $2$. 
To verify the signature only two multiplications and one square are needed.

Note that there is another signature scheme relying on the difficulty of finding square roots, the Rabin-Williams signature (cf. \cite{galbraith}), which avoids the forgery vulnerability. While that scheme requires the use of two primes respectively congruent to $3$ and $7$ modulo $8$, the two variants above do not need this condition. Moreover in the Rabin-Williams scheme, a message cannot be signed twice in two different ways, otherwise the factorization of $N$ might get exposed. In the above schemes, using a deterministic pad as above, allows different signatures of the same message.

For more on forgery and blindness on Rabin signatures, please refer also to \cite{DS}.

\section{Conclusions and Remarks}

Let us make here a few comments on the Rabin schemes in general, after having mainly dwelled on the deterministic aspects and identification problems.

In principle, the Rabin scheme is very efficient, because only one square is required for encryption;
 furthermore, it is provably as secure as factoring.
Nevertheless, it is well known \cite{buchmann,silverman} that it presents some drawbacks, mainly due to
 the four-to-one mapping, that may discourage its use to conceal the content of a message, namely: 
\begin{itemize}
  \item the root identification requires the delivery of additional information, which may increase computational costs; %not be easily
ù  \item many proposed root identification methods, based on the message semantics, have a probabilistic
   character and cannot be used in some circumstances;
  \item the delivery of two bits together with the encrypted message exposes the process to active attacks
   by maliciously  modifying these bits. For example, suppose an attacker $A$ sends an encrypted
    message to $B$ asking that the decrypted message be delivered to a third party $C$ (a friend of $A$).
    If in the encrypted message the bit that identifies the root among the two roots of the same parity had
    been deliberately changed, $A$ can get a root from $C$ that, combined with the original message, enables
    the Rabin public-key to be factored.
    Even Variant II is not immune to those kind of active attacks.
\end{itemize}

In conclusion, the Rabin scheme may suffer from some hindrance when used to conceal a message,
 whereas it seems effective when applied to generate an electronic signature or as a hash function. 
However, these observations do not exclude the practical use of the Rabin scheme (as is actually profitably done in some standardized protocols), when other properties, like integrity and authenticity,
    are to be taken care of, along with message secrecy, in a public-encryption protocol.

\section{Acknowledgments}
This work was partially done while the first author was Visiting Professor
with the University of Trento, funded by CIRM, and he would like to thank
 the Department of Mathematics for the friendly and fruitful atmosphere offered.
The third author has been supported by the Swiss National Science
Foundation under grant No. 132256. We would also like to thank Steven Galbraith for his comments on a preliminary version of the paper and for pointing out some references.

%\newpage 

% ******************************************************************

\end{document}